\documentclass[twoside]{article}
\usepackage{amsfonts}
\usepackage{alg}

\begin{document}

\title  {The polynomial identities for matrix functions}
\author{Egorychev~G.\\gegorych@mail.ru}
\maketitle



\textit{Abstract. }This is a short review of some recent results obtained by
the author [2 - 8]. These results are related the problem of obtaining
polynomial identities (computational formulas) for some matrix functions by
means of the known polarization theorem, including the case of
noncommutative variables and of determinant of space matrices.\bigskip 

\textbf{1. Introduction}

The classical polarization theorem [1; 2] about recovering of a polyadditive
symmetric function from its values on a diagonal plays a fundamental role in
the algebraic theory of polynomials. This theorem is often applied in the
ring theory, linear and polylinear algebra, theory of functions, the theory
of mixed volumes and in other areas of mathematics. We shall aplly it in the
following simple and general form.\bigskip

\textbf{Theorem 1 (}\textit{the polarization theorem}; \textit{Egorychev}{,
[2], 1980)}. \textit{Let} $H$\textit{\ be an additive commutative semigroup, 
}$H^{\left( n\right) }$\textit{\ }$=$\textit{\ }$H\times \ldots \times $%
\textit{\ }$H,$\textit{\ and }$\Phi $ \textit{be} \textit{an\ abelian group
with division by integers. The function $f(x)=f(x_{1},\ldots,x_{n})$: 
$H^{(n)}\rightarrow \Phi,$}\textit{\ is symmetric and polyadditive if and only if
for any }$\gamma $\textit{\ }$\in H$ \textit{it admits the following
representation} (\textit{the polarization formula}):\bigskip 
\[
f(x_{1},...,x_{n})=\{(-1)^{n}F(\gamma)
+\sum_{k=1}^{n}(-1)^{n-k}\sum_{1\leq j_{1}<...<j_{k}\leq n}F(\gamma
+\sum_{l=1}^{k}x_{j_{l}})\}/n!,
\]%
\textit{where}
\[
F(x):=f(x,\ldots ,x):H\rightarrow \Phi . 
\]%

Many well-known matrix functions (determinants, permanents, mixed
discriminants (mixed volumes in $\mathbb{R}^{n}$), resultants, etc.) are polyadditive functions of various types. Here
we shall show how in this case the polarization theorem can be used both for
calculation of these functions, and for defining them on more general
algebraic systems (with preservation of the basic properties).

Let $K$ be a commutative ring, $\mathbf{K}$ a noncommutative ring, $Q$ a
noncommutative ring with associative $n$-powers (one-monomial associativity)%
\textit{, }and let each ring be with division by integers.

\textbf{2. The first polynomial identity for permanents over a
commutative ring }

\bigskip
The concept of the permanent was first introduced independently and
practically simultaneously in the well known memoirs of J. Binet (1812) and
A. Cauchy (1812). It was at this time that J. Binet introduced the term
"permanent". The following classical definition of permanent over the
commutative ring (fields $\mathbb{R}$ and $\mathbb{C}$) is most often used.

\textbf{Definition 1}. \textit{If} $A=\left( a_{ij}\right) $ \textit{be an} $%
n\times n$ \textit{matrix} \textit{over a commutative ring}$\mathbf{,}$ 
\textit{then}
\[
per(A):=\sum_{\sigma \in S_{n}}a_{_{1}\sigma (1)}\times \ldots
\times a_{n\sigma(n)}, 
\]%
\textit{where }$S_{n}$\textit{\ be the set of all permutations }$\sigma
=\left( \sigma (1),\ldots,\sigma (n) \right) $%
\textit{\ of the set} $\{1,\ldots ,n\}.$

The author has for the first time obtained ({[3], 1979})\textit{\ } the
computation formula (polynomial identity) for permanents over a commutative
ring $K$ that contains up to $n!$ free variables. This result was obtained
by using the apparatus of the theory of power series. For an illustration we
give a short proof of the same result using the polarization formula.

\textbf{Theorem 2} (\textit{the first polynomial identity for permanents
over the commutative ring }$K$; \textit{Egorychev}{, [3], 1979}).
\[
{per(A)=\prod_{i=1}^{n}\gamma_{i}+\sum_{k=1}^{n}(-1)^{k}\sum_{1\leq j_{1}<...<j_{k}\leq n}\left( \prod_{i=1}^{n}
\left(\gamma_{i}-a_{ij_{1}}-\ldots -a_{ij_{k}}\right) \right),} 
\quad\left(2.1\right) 
\]%
\textit{where}{\ $\gamma _{1},\ldots ,\gamma _{n}$ \textit{are free elements
from} }$K$\textit{.}

\textbf{Proof}\textit{.}{\ Let }$A=${\ $(a_{1},\ldots ,a_{n})$, where $%
a_{1},\ldots ,a_{n}$ are its $n$-vectors-colunms. As the function }${per}${$%
(A)$ $=$\ }${per}${$(a_{1},\ldots ,a_{n})$ is a polyadditive and symmetric
function of $n$ variables $a_{1},\ldots ,a_{n},$ then by means of the
polarization formula for any $\gamma $ $=$ $\left( \gamma _{1},\ldots
,\gamma _{n}\right) $ }$\in ${\ }$K$$^{n}${\ we obtain
\[
per(A)=per(a_1,\ldots,a_n)=\{per (\gamma,\ldots,\gamma) + 
\]%
}%
\[
+\sum_{k=1}^{n}(-1)^{k}\sum_{1\leq j_{1}<...<j_{k}\leq n}per(\gamma -\sum_{l=1}^{k}a_{j_{l}},\ldots ,\gamma
-\sum_{l=1}^{k}a_{j_{l}})\}/n!= 
\]%
{\bigskip }(the equality $per (\gamma,\ldots,\gamma)/n!=\gamma_{1}\ldots \gamma_{n}$ for any elements $\gamma_{1},\ldots
,\gamma _{n}$ from $K)$
\[
{=\prod_{i=1}^{n}\gamma _{i}+\sum_{k=1}^{n}\left( -1\right) ^{k}\sum_{1\leq
j_{1}<...<j_{k}\leq n}\prod_{i=1}^{n}(\gamma_{i}-\sum_{l=1}^{k}a_{ij_{l}}).}%
\bullet 
\]%

\textbf{Remark 1.} \textit{Formula} (2.3) \textit{for special values of }$%
\gamma _{1},\ldots ,\gamma _{n}$ \textit{from} $K$ \textit{gives well-known
formulas of J. Riordan} (1966) \textit{and G. Wilf }(1968), \textit{which
are recognized today as the most economic for calculation of the permanent
over the commutative ring} $K$ (\textit{fields} {}$\mathbb{R}${}, {}$\mathbb{C}${}).

\textbf{3. Polynomial identity\ of a new type for determinants
over the commutative ring}

\bigskip 

Let $A=\left( a_{ij}\right) $ be an $n\times n$ matrix with elements from
the ring $K$. Let $S_{n}$ be the set of all permutations $\sigma =\left(
\sigma \left( 1\right) ,\ldots ,\sigma \left( n\right) \right) $ of the set $%
\{1,\ldots ,n\}.$ Let $\tau (\sigma )$ be a number of inversions in the
permutation $\sigma =$ $\left( \sigma \left( 1\right) ,\ldots ,\sigma \left(
n\right) \right) $ of the set $\{1,\ldots ,n\}.$ Let $S_{n}^{\left( e\right)
}$ and $S_{n}^{\left( o\right) }$ be the subsets of even and odd
permutations in $S_{n}$, respectively$.$ We call the sequence of elements of 
$a_{1\sigma \left( 1\right) },\ldots ,a_{n\sigma \left( n\right) }$ the
diagonal $l\left( \sigma \right) $ of matrix $A,$ and the sequence of
elements $a_{i_{1}\sigma \left( i_{1}\right) },\ldots ,a_{i_{k}\sigma \left(
i_{k}\right) },$ $1\leq i_{1}<\ldots <i_{k}\leq n$ the subdiagonal $l$ of
length $k$ of matrix $A$. We denote by $L_{k}^{\left( e\right) }\left(
L_{k}^{\left( o\right) }\right) $ the set of all subdiagonals of length $k$
of diagonals from $S_{n}^{\left( e\right) }\left( S_{n}^{\left( o\right)
}\right) ,$ $k=1,\ldots ,n.$The function $su\left( l\right) $ on the
subdiagonal (diagonal) $l=diag(a_{i_{1}\sigma \left( i_{1}\right) },\ldots
,a_{i_{k}\sigma \left( i_{k}\right) }),$ is the sum of its elements, i.e. $%
su\left( l\right) :=\sum_{s=1}^{k}a_{i_{s}\sigma \left( s\right) }.$

In [4] for the classical Cauchy determinant $\det \left( A\right) $ over a
commutative ring 

\[
\det \left( A\right):=\sum_{\sigma \in S_{n}}(-1)^{d(\sigma)}a_{1\sigma (1)}\times a_{2\sigma (2)}\ldots 
\times a_{n\sigma (n)} 
\]%
a series polynomial identities of a new type has been found. Among
them the following formula for a determinant which will be used in the
further statements.

\bigskip

\textbf{Theorem 3 (}\textit{polynomial identity\ for determinants over a
commutative ring }$K;$ \textit{Egorychev}\textbf{, }[4], 2012){.}

\[
\det \left( A\right) =\frac{1}{n!}\{[\sum_{l\in L_{n}^{(e)}}({%
\gamma +}su (l))^{n}-\sum_{l\in L_{n}^{(o)}}({%
\gamma +}su\left( l\right) )^{n}]- 
\]
\[
{-\ [\sum_{l\in L_{n-1}^{(e)}}({\gamma +}su (l))^{n}-\sum_{l\in L_{n-1}^{(o)}}({\gamma +}su (l))^{n}]\},}
\quad \quad\quad\,\, \left( 3.1\right) 
\]%
\textit{for any }$\gamma \in K.$\textit{\ In particular, for} $\gamma =0$%

\[
{\det \left( A\right) =\frac{1}{n!}\{[\sum_{l\in L_{n}^{\left( e\right)
}}su^{n}\left( l\right) -\sum_{l\in L_{n}^{\left( o\right) }}su^{n}\left(
l\right) ]-[\sum_{l\in L_{n-1}^{\left( e\right) }}su^{n}\left( l\right)
-\sum_{l\in L_{n-1}^{\left( o\right) }}su^{n}\left( l\right) ]\}.}
\,\,\left( 3.2\right) 
\]%

From (3.1) -- (3.2 ) the following corollary directly follows

\bigskip

\textbf{Corollary 1} ( \textit{Egorychev}\textbf{, }[4], 2012\textbf{).}%
\textit{If }$A=\left(à_{ij}\right) $\textit{\ is an }$n\times n$%
\textit{\ matrix over a commutative ring}$,$\textit{\ then for }$%
t=1,2,\ldots ,n-1$\textit{\ the following identities are valid}:
\[
\sum_{l\in L_{n}^{(e)}}su^{t}(l)-\sum_{l\in L_{n}^{(o)}}su^{t} (l)=\sum_{l\in L_{n-1}^{(e)}}su^{t}(l)-
\sum_{l\in L_{n-1}^{(o)}}su^{t}(l), 
\]%
\textit{and }$det(A)=0$\textit{\ iff}
\[
\sum_{l\in L_{n}^{(e)}}su^{n}(l)-\sum_{l\in L_{n}^{(o)}}su^{n} (l)=
\sum_{l\in L_{n-1}^{(e)}}su^{n}(l)-\sum_{l\in L_{n-1}^{(o)}}su^{n}(l).\quad \left( 3.3\right) 
\]%

\textit{In other words, vectors-rows }(\textit{vectors-columns}) 
\textit{of} \textit{an} $n\times n$ \textit{matrix} $A$ \textit{over a
commutative ring} $are$ $linearly$ $independent$ \textit{if and only if the
equality} (3.3) \textit{holds}.\bigskip 

\textbf{4. New formulas for permanents over noncommutative rings}\bigskip 

One of natural generalizations of a classic determinant on noncommutative
rings has been introduced in 2000 by A. Barvinok [9] who has successfully
applied it for calculation of permanents over a noncommutative ring. Similar
definition for \textit{the symmetrized\ permanent\ }$eper\left( A\right) $
has been given by the author in 2007.

\textbf{Definition 2 (}\textit{the symmetrized\ permanent }$eper\left(
A\right) $ \textit{over the noncommutative ring }$\mathbf{K}$; [5], 2007)$:$%

\[
eper\left( A\right) :=\sum_{\sigma \in S_{n}}Sym(a_{1\sigma (1)}
\times a_{2\sigma (2)}\ldots \times a_{n\sigma (n)
}), \quad \quad \left( 4.1\right) 
\]%
{\textit{where the symmetrization operator Sym for the product elements}} 
$x_1, x_2 ,\ldots, x_n$ $\in \mathbf{K}$ 
\textit{is the sum}

\[
Sym(x_{1}\times x_{2}\times\ldots \times x_{n}):=\frac{1}{n!}%
\sum_{\sigma \in S_{n}}x_{\sigma (1)}x_{\sigma (2)}\ldots x_{\sigma (n)}.\quad \quad \left( 4.2\right) 
\]%

In [4] the properties of the permanent $eper\left( A\right) $ over $\mathbf{K%
}$ have been studied, which coincide with many basic characteristic
properties of the classical $per\left( A\right) $ over $K$ [1]. Analogously
to [1], we further prove a polynomial identity for $eper\left( A\right) $
over $\mathbf{K}$ which generates a new family of identities for elements of
matrix $A$.

Let $A=\left(a_{ij}\right) $ be an $n\times n$\ matrix$,$ and let $\Phi_{r,s}(A)$ be the set of all own $r\times s$\ 
submatrices of the matrix $A$, 
${r,s}\in \{1,2,\ldots ,n\}.$ The function $su\left( A\right) $ of the matrix $%
A$ is the sum of its elements.

\textbf{Theorem 4 }\textit{\ If} $A=\left(a_{ij}\right) $ \textit{be an} $%
n\times n$ \textit{matrix} \textit{over the noncommutative ring }$\mathbf{K,}
$ then

\[
eper\left(A\right) =\frac{1}{n!}\sum_{r,s=1}^{n}(-1)^{r+s}
\sum\limits_{B\in \Phi_{r,s}(A)}(\delta +su\left(B\right))^{n}
\quad \left( 4.3\right) 
\]%
\textit{for any }$\delta \in \mathbf{K}.$%
\textit{\ In particular, for} $\delta =0$

\[
eper\left( A\right) =\frac{1}{n!}\sum_{r,s=1}^{n}(-1)^{r+s}
\prod\sum\limits_{B\in \Phi_{r,s}(A)}su^{n}\left( B\right). \quad \left( 4,4\right) 
\]%

From (4.3) -- (4.4 ) the following corollary directly follows.


\textbf{Corollary 2. }\textit{If }$A=\left(a_{ij}\right) $\textit{\
is an }$n\times n$\textit{\ matrix over a noncommutative ring} $\mathbf K$, \textit{
then for}$m=1,2,\ldots ,n-1$\textit{\ the following identities are valid}:%

\[
\sum_{r,s=1}^{n}(-1)^{r+s}
\sum\limits_{B\in \Phi_{r,s}(A)} su^{m}\left( B\right) =0,\, m=1,2,\ldots ,n-1, 
\]%
\textit{and }$edet(A)=0$\ 

\[
\sum_{r,s=1}^{n}(-1)^{r+s}\sum\limits_{B\in \Phi_{r,s}(A)}su^{n}\left(B\right) =0.
\]%

\textbf{Proof}. The $eper\left( A\right) $ is polyadditive and symmetric
function with respect to vector-columns of matrix $A.$ Applying\textit{\ }{%
the polarization formula to }$eper\left( A\right) $ we obtain

\[
eper\left( A\right) =\frac{1}{n!}\{\sum_{k=1}^{n}(-1)^{n-k}\sum_{1\leq
j_{1}<...<j_{k}\leq n}eper(\sum_{l=1}^{k}\alpha_{j_{ä }},\ldots
,\sum_{l=1}^{k}\alpha _{j_{l}})\}=\quad\left( 4.5\right) 
\]%
($eper\left( \alpha ,\alpha ,\ldots ,\alpha \right)
:=n!Sym\left( \alpha _{1}\times \alpha _{2}\times \ldots \times \alpha
_{n}\right) $ for any $n$-column $\alpha =\left( \alpha _{1},\alpha
_{2},\ldots ,\alpha _{n}\right) ^{T}\in $ $\mathbf{K}^{n}$)

\[
=\{\sum_{k=1}^{n}(-1)^{n-k}\sum_{1\leq j_{1}<...<j_{k}\leq n}
Sym(\prod\limits_{i=1}^{n}(\sum_{l=1}^{k}a_{ij_{l}})\}. 
\]%
Applying\textit{\ }{the polarization formula to }$%
Sym(\prod\limits_{i=1}^{n}(\sum_{l=1}^{k}a_{ij_{l}})$, we obtain 
\[
Sym\{\prod\limits_{i=1}^{n}(\sum_{l=1}^{k}a_{ij_{l}})\}=\frac{1}{n!}%
\{(-1)^{n}Sym \left(\delta^{n}\right) + 
\]
\[
+\sum_{s=1}^{n}(-1)^{n-s}\sum_{1\leq i_{1}<...<i_{s}\leq n}Sym((\delta
+\sum_{t=1}^{s}\sum_{l=1}^{k}a_{ij_{l}})^{n})\}. 
\]%
($Sym\left(x^{n}\right) =x^{n}$ for any $x\in \mathbf{K}$)

\[
=\frac{1}{n!}\{(-1)^{n}\delta^{n}+\sum_{s=1}^{n}(-1)^{n-s}\sum_{1\leq i_{1}<...<i_{s}\leq n}(\delta
+\sum_{t=1}^{s}\sum_{l=1}^{k}a_{ij_{l}})^{n}\}.\quad \left(4.6\right) 
\]%

Now it is enough to substitute in $\left( 4.5\right) $ the
expression $\left( 4.6\right) $ for $Sym(\prod\limits_{i=1}^{n}(%
\sum_{l=1}^{k}a_{ij_{l}})$, and to notice, that the coefficient at $\delta^n $ 
is equal to zero.$\bullet $

\bigskip 
\textbf{Remark 2. }\textit{In the proof of Theorem 4 we have chosen the
values of free elements from }$\mathbf{K}$\textit{\ to be equal to zero, or
to the element }$\delta $\textit{. Otherwise, we obtain a family of formulas
the more complicated the more free elements are in }$\mathbf{K}$\textit{.
Each of these formulas can be taken as the definition for }$eper\left(
A\right) $\textit{\ and generates a new family of identities for elements of
the matrix} $A$.

\bigskip 

\textbf{5. \ Determinants of space matrices over a commutative ring}\bigskip 

\textbf{Definition 3. }\textit{Let} $\widetilde{A}=(a_{ij}^{\left( k\right)
})$ be an $n\times n\times n$ \textit{matrix over the commutative ring }$K,$%
\textit{\ and let\ }$\widetilde{A}=\left(A_{1}, A_{2},\ldots, A_{n}\right) ,$\textit{\ where }
$A_{k}=(a_{ij}^{\left(k\right) }),$\textit{\ }$k=1,2,\ldots ,n,$\textit{\ are }$n\times n$\textit{%
\ sections of matrix }$\widetilde{A}$\textit{\ at index }$k.$ \textit{\ In
turn let}$\ A_{k}=(\alpha_{1}^{\left( k\right) },\alpha
_{2}^{\left( k\right) },\ldots ,\alpha _{n}^{\left( k\right) })$, where $%
\alpha _{j}^{\left( k\right) }=(a_{1j}^{\left( k\right) },a_{2j}^{\left(
k\right) },\ldots ,a_{nj}^{\left( k\right) })^{T},$ $j=1,2,\ldots ,n,$ 
\textit{are }$n\times n$\textit{\ vectors-columns of matrix }$A_{ê }$. 
\textit{Let's define matrix function} $Det_{p}(\widetilde{A})$ 
\textit{as follows} (\textit{comp. }[9]\textit{, }pp. 11--25):

\[
Det_{p}(\widetilde{A}):=\sum_{\sigma \in S_{n}}\left( -1\right) ^{d\left(
\sigma \right) }per(\alpha _{\sigma \left( 1\right) }^{\left( 1\right)
},\alpha _{\sigma \left( 2\right) }^{\left( 2\right) },\ldots ,\alpha
_{\sigma \left( n\right) }^{\left( n\right) }). 
\]%

\textbf{Theorem 5.} \textit{The following formula is valid }[8]:

\[
Det_{p}(\widetilde{A})=\{[\sum_{\sigma \in S_{n}^{\left( e\right)
}}\prod\limits_{t=1}^{n}(\sum_{i=1}^{n}a_{t\sigma \left( i\right) }^{\left(
i\right) })-\sum_{\sigma \in S_{n}^{\left( o\right)
}}\prod\limits_{t=1}^{n}(\sum_{i=1}^{n}a_{t\sigma \left( i\right) }^{\left(
i\right) })]- 
\]

\[
-[\sum_{\sigma \in S_{n}^{\left( e\right)
}}(\sum_{s=1}^{n}\prod\limits_{t=1}^{n}(-a_{t\sigma \left( s\right)
}^{\left( s\right) }+\sum_{i=1}^{n}a_{t\sigma \left( i\right) }^{\left(
i\right) }))-\sum_{\sigma \in S_{n}^{\left( o\right)
}}(\sum_{s=1}^{n}\prod\limits_{t=1}^{n}(-a_{t\sigma \left( s\right)
}^{\left( s\right) }+\sum_{i=1}^{n}a_{t\sigma \left( i\right) }^{\left(
i\right) }))]. 
\]%

\textbf{Proof.}\textit{\ }Applying\textit{\ }{the polarization formula to a
general member of the determinant sum for }$Det_{p}(\widetilde{A}),${\ we
obtain}

\[
per(\alpha _{\sigma \left( 1\right) }^{\left( 1\right) },\alpha _{\sigma
\left( 2\right) }^{\left( 2\right) },\ldots ,\alpha _{\sigma \left( n\right)
}^{\left( n\right) })= 
\]

\[
\frac{1}{n!}\{[\sum_{l\in L_{n}^{\left( e\right) }}per(\sum_{i=1}^{n}\alpha
_{\sigma \left( i\right) }^{\left( i\right) },\ldots ,\sum_{i=1}^{n}\alpha
_{\sigma \left( i\right) }^{\left( i\right) })-\sum_{l\in L_{n}^{\left(
o\right) }}per(\sum_{i=1}^{n}\alpha _{\sigma \left( i\right) }^{\left(
i\right) },\ldots ,\sum_{i=1}^{n}\alpha _{\sigma \left( i\right) }^{\left(
i\right) })], 
\]%
and
\[
Det_{p}(\widetilde{A})=\frac{1}{n!}\{[\sum_{l\in L_{n}^{\left( e\right)
}}per(\sum_{i=1}^{n}\alpha _{\sigma \left( i\right) }^{\left( i\right)
},\ldots ,\sum_{i=1}^{n}\alpha _{\sigma \left( i\right) }^{\left( i\right)
})-\sum_{l\in L_{n}^{\left( o\right) }}per(\sum_{i=1}^{n}\alpha _{\sigma
\left( i\right) }^{\left( i\right) },\ldots ,\sum_{i=1}^{n}\alpha _{\sigma
\left( i\right) }^{\left( i\right) })]- 
\]
\[
-[\sum_{s=1}^{n}\sum_{l\in L_{n}^{\left( e\right) }}per(-\alpha _{\sigma
\left( s\right) }^{\left( s\right) }+\sum_{i=1}^{n}\alpha _{\sigma \left(
i\right) }^{\left( i\right) },\ldots ,-\alpha _{\sigma \left( s\right)
}^{\left( s\right) }+\sum_{i=1}^{n}\alpha _{\sigma \left( i\right) }^{\left(
i\right) })- 
\]
\[
-\sum_{s=1}^{n}\sum_{l\in L_{n}^{\left( o\right) }}per(-\alpha _{\sigma
\left( s\right) }^{\left( s\right) }+\sum_{i=1}^{n}\alpha _{\sigma \left(
i\right) }^{\left( i\right) },\ldots ,-\alpha _{\sigma \left( s\right)
}^{\left( s\right) }+\sum_{i=1}^{n}\alpha _{\sigma \left( i\right) }^{\left(
i\right) })]\}= 
\]%
 (the equality $per\left( \gamma ,\ldots ,\gamma \right) /n!=\gamma
_{1}\ldots \gamma _{n}$ for any elements $\gamma _{1},\ldots ,\gamma _{n}$
from $K)$
\[
=\{[\sum_{\sigma \in S_{n}^{\left( e\right)
}}\prod\limits_{t=1}^{n}(\sum_{i=1}^{n}a_{t\sigma \left( i\right) }^{\left(
i\right) })-\sum_{\sigma \in S_{n}^{\left( o\right)
}}\prod\limits_{t=1}^{n}(\sum_{i=1}^{n}a_{t\sigma \left( i\right) }^{\left(
i\right) })]- 
\]
\[
-[\sum_{\sigma \in S_{n}^{\left( e\right)
}}(\sum_{s=1}^{n}\prod\limits_{t=1}^{n}(-a_{t\sigma \left( s\right)
}^{\left( s\right) }+\sum_{i=1}^{n}a_{t\sigma \left( i\right) }^{\left(
i\right) }))-\sum_{\sigma \in S_{n}^{\left( o\right)
}}(\sum_{s=1}^{n}\prod\limits_{t=1}^{n}(-a_{t\sigma \left( s\right)
}^{\left( s\right) }+\sum_{i=1}^{n}a_{t\sigma \left( i\right) }^{\left(
i\right) }))].\bullet 
\]

\textbf{Remark 3. }\textit{The} \textit{formulas of theorems 3--5\ and its
corollaries employ }(\textit{besides divisions by }$n!$)\textit{\ \ only
operations }$+$\textit{,}$-$\textit{,\ and raising to power }$n.${\large \ }%
\textit{Therefore} \textit{these} \textit{formulas} \textit{can be used} 
\textit{as a definition} \textit{for }$e$-\textit{determinants }($e$-\textit{%
permanents}) \textit{over a noncommutative ring }$Q$ \textit{with
associative }$n$-\textit{powers }(\textit{with preservation of the basic
properties} \textit{of these matrix functiîns}).\textbf{\ }

Obtaining of similar results is of much interest for the Shur functions,
mixed discriminants and many other matrix functions of planar and space
matrices. The direct application of our results can be found in the theory
of permanents [7], tensor algebra and its applications, the theory of $n$%
-Lie algebras [11; 12], the theory of noncommutative skew fields [13],
differential geometry [14], and others.

The author is grateful to his colleagues S. G. Kolesnikov, V. M. Kopytov,
V.P. Krivokolesko, V. M. Levchuk, A.P. Pozhidaev, A,V. Schuplev for
discussion of the basic results of this work and a number of useful
remarks.

\begin{center}
\textbf{References}
\end{center}

1. Cartan H. Elementary theory of analytic functions of one or several
complex variables. N. Y., Dover Publ., 1995. 300 p.

2. Egorychev G.P. New formulas for the permanent. Soviet Math. Dokl. 1980,
vol. 265, no 4, pp. 784--787 (in Russian).

3. Egorychev G.P. A polynomial identity for the permanent. Math. Zametki.,
1979, vol. 26, no 6, pp. 961--964. (in Russian)

https://doi.org/10.1007/BF01142089

4. Egorychev G.P. New polynomial identities for determinants over commutative 
rings. Izv. Irkutsk. Gos. Univ. Ser. Mat., 2012, vol. 5, no 4, pp. 16--20
(in Russian)

5. Egorychev G.P., Egorycheva Z.V. New polynomial identity for computing 
permanents. Sibirian Federal University. Dep. in VINITI
13.06.2007, no 627-B2007, pp. 1--32 (in Russian).

6. Egorychev G.P. A new family of polynomial identities for computing
determinats. Dokl. Acad. Nauk, 2013, vol. 452, no 1, pp. 14--16. (in Russian)

7. Egorychev G.P. Discrete mathematics. Permanents. Krasnoyarsk, Siberian
Federal University Publ., 2007. 272 p.

8. Egorychev G.P. The polynomial identities for matrix functions. Izv.
Irkutsk. Gos. Univ. Ser. Mat., 2017, vol. 21, no 4, pp. 77--88 (in Russian).

9. Barvinok A. New permanent estimators via non-commutative determinants. 
arXiv preprint math/0007153, 2000, arXiv:math/0007153, 2000, pp. 1--13.

10. Sokolov N.P. Spañe matrices and its applications. Moscow,
PhyzMathLit Publ., 1960. 300 p.

11. Filippov V.T. On the $n$-Lie algebras of Jacobian algebras . Sibirsk.
Mat. Zh., 1998, vol. 39, no 3, pp. 660--669.(in Russian)

12  Pozhidaev A.P. A simple factor-algebras and subalgebras of Jacobian
algebras, Sibirsk.Mat. Zh., \textbf{39}(1998), No. 3, 593--599 (Russian).

13.\ Kolesnikov P.S. About different definitions of algebraically closed
skew fields, Algebra and Logic, \textbf{40}(2001), No. 4, 396--414 (Russian).

14. Sabinin L.V. Methods of Nonassociative Algebra in Differential Geometry.
In Supplement to Russian translation of S.K. Kobayashi and K. Nomizu
Foundations of Differential Geometry,1, Moscow, Nauka, 1981 (Russian)

\end{document}